\documentclass[ reqno]{amsart}
\usepackage{amsmath, amssymb, amsthm, geometry, enumerate, graphicx, amsfonts, hyperref, mathrsfs}
\hypersetup{colorlinks=true,linkcolor=red, anchorcolor=blue, citecolor=blue, urlcolor=red, filecolor=magenta, pdftoolbar=true}
\theoremstyle{plain}
\newtheorem{thm}{\bf Theorem}[section]
\newtheorem{prop}[thm]{\bf Proposition}

\theoremstyle{remark}
\newtheorem{defn}[thm]{\bf Def{}inition}
\newtheorem{rem}[thm]{\bf Remark}
\newtheorem{exa}[thm]{\bf Example}
\numberwithin{equation}{section}
\usepackage{enumerate}
\usepackage{graphicx}
\allowdisplaybreaks
\begin{document}
\baselineskip08pt
\title [Weaving $K$-Frames in Hilbert Spaces] {Weaving $K$-Frames in Hilbert Spaces}

\author[Deepshikha]{Deepshikha}
\address{Deepshikha, Department of Mathematics,
University of Delhi, Delhi-110007, India.}
\email{dpmmehra@gmail.com}

\author[Lalit K. Vashisht]{Lalit  K. Vashisht}
\address{Lalit  K. Vashisht, Department of Mathematics,
University of Delhi, Delhi-110007, India.}
\email{lalitkvashisht@gmail.com}

\begin{abstract}
  G$\textnormal{\v{a}}$vruta introduced   $K$-frames for Hilbert spaces  to study atomic systems with respect to a bounded linear operator. There are many differences between $K$-frames and standard
frames, so we study weaving properties of $K$-frames.  Two frames $\{\phi_{i}\}_{i \in I}$  and $\{\psi_{i}\}_{i \in I}$ for a separable  Hilbert space $\mathcal{H}$ are  woven if there are positive constants $A \leq B$ such that for every subset $\sigma \subset I$, the family $\{\phi_{i}\}_{i \in \sigma} \cup \{\psi_{i}\}_{i \in \sigma^{c}}$  is a frame for $\mathcal{H}$ with frame bounds $A,  B$.  In this paper,  we present necessary and sufficient conditions for weaving  $K$-frames in   Hilbert spaces.  It is shown that woven $K$-frames and weakly woven $K$-frames are equivalent.   Finally, sufficient conditions for Paley-Wiener type perturbation of  weaving $K$-frames are given.
\end{abstract}


\keywords{  Frames; $K$-frames;   weaving frames; local atoms; perturbation.\\
 The research of Deepshikha is supported by the Council of Scientific $\&$ Industrial Research (CSIR) (Grant No.: 09/045(1352)/2014-
EMR-I), India. Lalit was   supported by R $\&$ D Doctoral Research Programme, University of Delhi
(Grant No.: RC/2015/9677).\\
$^*$Corresponding author-Lalit K. Vashisht}

\maketitle
 \baselineskip12pt

\section{Introduction and Preliminaries}
Let $\mathcal{H}$ be a  complex separable   Hilbert space with an inner product $\langle ., . \rangle$. A countable   sequence $\{f_k\}_{k=1}^{\infty} \subset \mathcal{H}$   is a \emph{frame}  for  $\mathcal{H}$  if there exist  positive scalars $A_o \leq B_o < \infty$ such that
\begin{align}\label{1.1}
A_o \|f\|^2\leq  \sum_{k = 1}^{\infty} |\langle f, f_k\rangle|^2 \leq B_o \|f\|^2
\  \text{for all} \  f \in \mathcal{H}.
\end{align}
The scalars $A_o$ and $B_o$ are called \emph{lower} and \emph{upper frame bounds}, respectively.  The sequence  $\{f_k\}_{k=1}^{\infty}$ is called  a \emph{Bessel sequence}  with \emph{Bessel bound} $B_o$ if  the upper inequality in \eqref{1.1} holds for all $f \in \mathcal{H}$.

Following three  operators  are associated with a frame $\{f_k\}_{k=1}^{\infty}$  for $\mathcal{H}$:
\begin{align*}
& \text{\emph{pre-frame  operator}} \ \    T:\ell^2(\mathbb{N})\rightarrow \mathcal{H}, \quad  \ T \{c_k\}_{k=1}^{\infty} =\sum\limits_{k = 1}^{\infty} c_k f_k, \  \{c_k\}_{k=1}^{\infty} \in \ell^2(\mathbb{N}),\\
& \text{\emph{analysis operator \ $($adjoint of $T)$}} \ \   T^{*}:\mathcal{H} \rightarrow
\ell^2(\mathbb{N}),  \quad   T^{*} f =\{\langle f, f_k\rangle\}_{k=1}^{\infty},  f \in \mathcal{H},\\
& \text{\emph{frame operator}} \quad S = T T^{*}: \mathcal{H}\rightarrow
\mathcal{H}, \quad   S f = \sum\limits_{k = 1}^{\infty}\langle f, f_k\rangle f_k,
f \in  \mathcal{H}.
\end{align*}
The frame operator $S$  is a  bounded, linear and invertible operator on $\mathcal{H}$. This gives the \emph{reconstruction} of each vector  $f \in \mathcal{H}$,
\begin{align*}
&f = SS^{-1}f =\sum\limits_{k = 1}^{\infty} \langle S^{-1}f, f_k \rangle f_k.
\end{align*}
Thus,  a frame for $\mathcal{H}$ allows each vector
in $\mathcal{H}$  to be written as a linear combination of the elements in the frame, but the linear independence between the frame
elements is not required. Concerning the evolution of the notion of Hilbert frames and their applications in different directions in science and engineering, it is necessary to mention the nobel books by Casazza and Kutyniok  \cite{CK},  Christensen \cite{OC}, Gr\textnormal{\"{o}}chenig \cite{Gro}; and  beautiful research tutorials  by Casazza \cite{C} and Casazza and  Lynch \cite{CL1}.

Next we give some basic notations.  The family of all bounded linear operators from a Banach space $\mathcal{X}$ into a Banach space $\mathcal{Y}$ is denoted by $\mathcal{B}(\mathcal{X}, \mathcal{Y})$. If $\mathcal{X} = \mathcal{Y}$, then we write $\mathcal{B}(\mathcal{X}, \mathcal{Y})= \mathcal{B}(\mathcal{X})$. The range and the
kernel of  $K \in \mathcal{B}(\mathcal{X}, \mathcal{Y})$ are denoted by $R(K)$ and $N(K)$, respectively. The pseudo-inverse of $K \in \mathcal{B}(\mathcal{H})$ is denoted by $K^{\dagger}$. Note that $K K^{\dagger} f = f$ for all $f \in R(K)$. Throughout the paper $R(K)$ is closed. By $\mathbb{N}$ we denote the set of all positive integers. The canonical orthonormal basis for $\ell^2(\mathbb{N})$  is the sequence $\{e_n\}_{n=1}^{\infty}$, where $e_n = \{0,0,0, \dots,\underbrace{1}_{nth \ place},0,0,0, \dots\}$ for all  $n \in \mathbb{N}$. For a sequence of vectors $\{f_k\}_{k \in I} \subset \mathcal{H}$, the closure of the span$\{f_k\}_{k \in I}$ is denoted by $[f_k]_{k \in I}$.\\

\vspace{5pt}
The following key-theorem can be found in \cite{Do}
\begin{thm}\cite{Do}\label{th1.3}
Let $\mathcal{H}_1$ and $\mathcal{H}_2$ be separable Hilbert spaces and let $L_1\in \mathcal{B}(\mathcal{H}_1,\mathcal{H})$, \break $L_2\in \mathcal{B}(\mathcal{H}_2,\mathcal{H})$. The following statements are equivalent:
\begin{enumerate}[$(i)$]
\item $R(L_1)\subset R(L_2)$
\item $L_1L_1^*\leq\lambda^2L_2L_2^*$ for some $\lambda\geq0$ and
\item there exists a  $C\in\mathcal{B}(\mathcal{H}_1,\mathcal{H}_2)$ such that $L_1=L_2C$ and
 $\|C\|^2=\inf\{\mu:L_1L_1^*\leq\mu L_2L_2^*\}$.
\end{enumerate}
\end{thm}

\subsection{$K$-frames in Hilbert spaces}
Feichtinger and Werther  \cite{FW}  introduced a family of analysis and synthesis systems with frame-like properties  for closed subspaces of  $\mathcal{H}$ and call it an \emph{atomic system} (or \emph{local atoms}). The motivation for the atomic system is based on examples arising in sampling theory, see \cite{FG}. One of the important properties of the atomic system is that it can generate a proper subspace even though they do not belong to them.

 \begin{defn}\cite{FW}
 Let $\mathcal{H}_o$ be a closed subspace of  $\mathcal{H}$. A sequence $\{f_k\}_{k=1}^{\infty} \subset
\mathcal{H}$ is called a \emph{family of local atoms} (or \emph{atomic system}) for
$\mathcal{H}_o$, if
\begin{enumerate}[$(i)$]
 \item  $\{f_k\}_{k=1}^{\infty}$ is Bessel sequence in $\mathcal{H}$,
\item there exists a sequence of linear functionals $\{c_k\}_{k=1}^{\infty}$ and  a real number  $C > 0$ such that
\begin{align}
 \sum_{k=1}^{\infty} |c_k(f)|^2  \leq C \|f\|^2 \ \text{for  all} \  f \in \mathcal{H}_o, \notag
\end{align}
\item\label{eq1.2}
$f = \sum_{k=1}^{\infty} c_k(f) f_k \ \text{for  all} \ f \in \mathcal{H}_o.$
\end{enumerate}
\end{defn}
\begin{rem}
Note that the linear functionals $\{c_k\}_{k=1}^{\infty}$ need only to be defined on the subspace $\mathcal{H}_o$. We say that $c_k$ are \emph{associated functionals} of the local atoms $\{f_k\}$. The constant $C$ is called the  \emph{atomic bound}. Furthermore, the partial sum $\sum_{k=1}^{N}c_k(f) f_k$ of the series in \eqref{eq1.2} can be converges to $f$ from ``outside'' of $\mathcal{H}_o$. The family of local atoms for $\mathcal{H}_o$ is a so-called
\emph{pseudo-frame} as proposed  by  Li and Ogawa in \cite{LO}.
\end{rem}

 G$\textnormal{\v{a}}$vruta in \cite{LG}  introduced and studied $K$-frames
in Hilbert spaces to study atomic systems with respect to a bounded
linear operator $K$ on Hilbert spaces.
\begin{defn}\cite{LG}
Let $K \in \mathcal{B}(\mathcal{H})$.  A sequence $\{f_k\}_{k=1}^{\infty} \subset
\mathcal{H}$ is called a \emph{$K$-frame} for $\mathcal{H}$, if there exist
 constants $A, B
> 0$ such that
\begin{align}
A\|K^*f\|^2\leq \sum_{k=1}^{\infty} |\langle f, f_k\rangle|^2 \leq
B\|f\|^2 \ \text{for all} \ f \in \mathcal{H}.
\end{align}
\end{defn}

 The numbers $A$ and $B$ are called \emph{lower} and \emph{upper} \emph{$K$-frame bounds}, respectively. If $I$ is the identity operator on $\mathcal{H}$, then $K$-frames are the standard frames.  $K$-frames are more general than standard frames in the sense that the lower frame bound only holds for the elements in the range of $K^*$. Since a $K$-frame $\{f_k\}_{k=1}^{\infty}$ for $\mathcal{H}$ is a Bessel sequence, we can define the pre-frame operator, analysis operator and  frame operator associated with  $\{f_k\}_{k=1}^{\infty}$. The frame operator of a $K$-frame is not
invertible on $\mathcal{H}$ in general, but  it is invertible on a subspace $R(K)$, where the range $R(K) \subset \mathcal{H}$ is closed. Furthermore, there are many differences between $K$-frames and standard frames, see \cite{LG}.  G$\textnormal{\v{a}}$vruta  \cite{LG}  characterize $K$-frames in
Hilbert spaces by using bounded linear operators. In \cite{DV3, LG2, LG3, XZG} some new results about $K$-frames  were obtained.

\subsection{Background on weaving frames}
The concept  of \emph{``weaving frames''} for complex separable Hilbert spaces introduced  by  Bemrose,  Casazza,  Gr\"{o}chenig,   Lammers  and  Lynch  in \cite{BCGLL}.
For a fixed $m \in \mathbb{N}$, we write $[m] = \{1,2, \dots, m\} \quad \text{ and } \ [m]^{c} = \mathbb{N} \setminus [m] = \{m+1, m+2, \dots\}$.

\begin{defn}\cite{BCGLL}
Let $I$ be a countable  indexing set. A family of frames $\Big\{\{\phi_{ij}\}_{j\in I}: i \in [m]\Big\}$  for  $\mathcal{H}$ is said to
be \emph{woven}, if there are universal constants $A$ and $B$ such  that for every partition $\{\sigma_i\}_{i \in [m]}$ of $I$,  the family $\bigcup_{i\in[m]}\{\phi_{ij}\}_{j\in\sigma_i}$  is a frame for $\mathcal{H}$ with  frame bounds $A$
and $B$.
\end{defn}
This new   notion  of weaving frames is  motivated by a problem in distributed signal processing. Weaving frames has potential applications in wireless sensor networks that require distributed processing under different frames, as well as pre-processing of signals using Gabor frames. Bemrose et al. \cite{BCGLL} gave a characterization of  weaving frames (that does not require universal frame bounds) and weaving Riesz bases. They  proved a  geometric characterization of woven Riesz bases in terms of distance between subspaces of a Hilbert space. Sufficient
conditions for weaving frames by means of perturbation theory and diagonal dominance can be found in \cite{BCGLL}.
The   fundamental properties of weaving frames  reviewed by Casazza and Lynch in \cite{CL}.  Casazza and Lynch  \cite{CL} proved some basic properties in the theory of weaving frames. They proved that an invertible operator applied to woven frames leaves them woven.  Casazza and Lynch \cite{CL} considered a ``weaving equivalent'' of an unconditional basis for weaving Riesz basis.   Casazza, Freeman and Lynch \cite{CFL} extended the concept of weaving Hilbert space frames to the Banach space setting. They introduced and studied  \emph{weaving Schauder frames} in Banach spaces. It is  proved in \cite{CFL}  that  for any two approximate Schauder frames for a Banach space $\mathcal{X}$, every weaving is an approximate Schauder frame if and only if there is a uniform constant $C \geq 1$ such that every weaving is a $C$-approximate Schauder frame. Some perturbation theorems for woven approximate Schauder frames can be found in  \cite{CFL}.   Deepshikha and Vashisht  studied weaving properties of an infinite family of frames in separable Hilbert spaces  in \cite{DV}. In \cite{VD, VD1}, authors introduced and studied  weaving frames  with respect to  measure spaces. Weaving properties of generalized frames and fusion frames can be found in \cite{DVV, DVV17, VGDD}. In  \cite{MoMF},  D\"{o}rfler and  Faulhuber studied weaving Gabor frames in $L^2(\mathbb{R})$. They gave a sufficient criteria for
a family of multi-window Gabor frames to be woven. A family of  localization operators related to weaving Gabor frames  is also discussed in \cite{MoMF}.
\subsection{Outline of the paper}
In this paper  we give  necessary and sufficient conditions for weaving  $K$-frames in   Hilbert spaces. A characterization of weaving  $K$-frames in terms of a bounded linear operator is given, see Theorem \ref{th2}. In Theorem \ref{th2.7}, a  sufficient condition for  $K$-frames  not to be weakly woven  in terms of  lower $K$-frame bounds is given. Theorem \ref{th2.6n} shows that woven $K$-frames and weakly woven $K$-frames are equivalent.  A characterization of weaving $K$-frames in terms of action of operators on $K$-frames is presented, see Proposition \ref{pr2}. To be precise, it is shown that an operator $U \in \mathcal{B}(\mathcal{H})$ applied to woven  $K$-frames leaves them $UK$-woven. A Paley-Wiener type    perturbation result for   weaving $K$-frames is  given in Theorem \ref{th8}. Several examples and counter-examples are given to illustrate the results.

\section{Main Results}
We begin with the definition.
\begin{defn}
 A family of  $K$-frames $\Big\{\{\phi_{ij}\}_{j=1}^\infty:i\in[m]\Big\}$ for  $\mathcal{H}$  is said to be    \emph{K-woven} if   there exist universal positive  constants $A$ and $B$ such that for any partition $\{\sigma_i\}_{i\in [m]}$ of $\mathbb{N}$, the family $\bigcup_{i\in[m]}\{\phi_{ij}\}_{j\in\sigma_i}$ is a $K$-frame for $\mathcal{H}$ with lower and upper  $K$-frame bounds  $A$ and $B$, respectively. Each family  $\bigcup_{i\in[m]}\{\phi_{ij}\}_{j\in\sigma_i}$ is  called a weaving.
\end{defn}

\begin{defn}
A family of  $K$-frames $\Big\{\{\phi_{ij}\}_{j=1}^\infty:i\in[m]\Big\}$ for  $\mathcal{H}$  is said to be    \emph{weakly K-woven}  if    for any partition $\{\sigma_i\}_{i\in [m]}$ of $\mathbb{N}$, the family $\bigcup_{i\in[m]}\{\phi_{ij}\}_{j\in\sigma_i}$ is a $K$-frame for $\mathcal{H}$.
\end{defn}

In \cite[Theorem 4.5]{BCGLL}, Bemrose et al. proved that if each weaving is a frame for $\mathcal{H}$,  then
there exist uniform frame bounds $A_o \leq B_o < \infty$ that simultaneously work for
all weavings. To be precise, weakly woven is equivalent to the frames being woven. We extend this result for  $K$-frame, see Theorem \ref{th2.6n}. Our approach is different from the technique given in \cite[Theorem 4.5]{BCGLL}, see Remark \ref{rem9ext1} for the details. In the direction of upper universal bounds for ordinary frames, Bemrose et al.  showed in  \cite[Proposition 3.1]{BCGLL} that  every weaving of ordinary frames automatically has a
universal upper frame bound. Using a technique given in  \cite[Proposition 3.1]{BCGLL}, we have the following result for universal upper frame bounds for $K$-frames.
\begin{prop}\label{pro1}
For  each $i \in [m]$, let   $\{\phi_{ij}\}_{j=1}^\infty$  be  a $K$-frame for $\mathcal{H}$ with $K$-frame bounds $A_i$ and $B_i$. Then, for any partition $\{\sigma_i\} _{i\in [m]}$ of $\mathbb{N}$, the family $\bigcup_{i\in[m]}\{\phi_{ij}\}_{j\in\sigma_i}$ is a Bessel sequence with Bessel  bound $\sum_{i\in[m]} B_i$. That is,  $\sum_{i\in[m]} B_i$ is one of the choice for an universal upper $K$-frame bound.
\end{prop}
\proof
Let $\{\sigma_i\} _{i\in[m]}$ be any partition of $\mathbb{N}$. Then, for any $f\in\mathcal{H}$, we have
\begin{align*}
\sum\limits_{i\in[m]} \sum\limits_{j\in\sigma_i}|\langle f,\phi_{ij}\rangle|^2 \leq\sum\limits_{i\in[m]} \sum\limits_{j\in\mathbb{N}}|\langle f,\phi_{ij}\rangle|^2
\leq \Big(\sum\limits_{i\in[m]} B_i \Big)\|f\|^2.
\end{align*}
This gives the required universal upper  $K$-frame bound for the family  $\bigcup_{i\in[m]}\{\phi_{ij}\}_{j\in\sigma_i}$.
\endproof

The following theorem gives a necessary and sufficient condition for  weaving $K$-frames in terms of an operator.
\begin{thm}\label{th2}
For each $i\in [m]$, suppose $\{\phi_{ij}\}_{j=1}^\infty$ is a $K$-frame for $\mathcal{H}$ with  bounds $A_i$ and $B_i$. The following conditions are equivalent.
\begin{enumerate}[$(i)$]
 \item The family $\Big\{\{\phi_{ij}\}_{j=1}^\infty:i\in[m]\Big\}$   is $K$-woven.
 \item  There exists $A>0$ such that for any partition  $\sigma = \{\sigma_i\}_{i\in[m]}$ of $\mathbb{N}$ there exists a bounded linear operator $M_{\sigma}:\ell^2(\mathbb{N})\rightarrow\mathcal{H}$ such that
\begin{align*}
M_{\sigma}(e_j)= &\begin{cases}
\phi_{1j}, j\in\sigma_1,\\
\phi_{2j}, j\in\sigma_2,\\
\vdots\\
\phi_{mj}, j\in\sigma_m,
\end{cases}
\end{align*}
\end{enumerate}
and $AKK^{*}\leq M_{\sigma}M_{\sigma}^{*}$, where $\{e_j\}_{j=1}^{\infty}$ is the  canonical orthonormal  basis for $\ell^2(\mathbb{N})$.
\end{thm}
\proof
$(i) \Rightarrow (ii):$ Suppose $A$ is an  universal lower $K$-frame bound for  the family \break $\Big\{\{\phi_{ij}\}_{j=1}^\infty:i\in[m]\Big\}$.
For any partition  $\sigma = \{\sigma_i\}_{i\in[m]}$ of $\mathbb{N}$, let $T_{\sigma}$ be the pre-frame operator associated with the Bessel sequence $\bigcup_{i\in[m]}\{\phi_{ij}\}_{j\in\sigma_i}$.

Choose $M_{\sigma}=T_{\sigma}$.
Then, $M_{\sigma}(e_j)=T_{\sigma}(e_j)=\phi_{ij}$ for all $j \in \sigma_i$ $(i \in [m])$.\\
Next, we compute
\begin{align*}
A\langle KK^{*}f,f\rangle &= A\|K^{*}f\|^2\\
&\leq\sum\limits_{i\in[m]}\sum\limits_{j\in \sigma_{i}}|\langle f,\phi_{ij}\rangle|^2\\
&=\sum\limits_{j\in\mathbb{N}}|\langle f,M_{\sigma}(e_j)\rangle|^2\\
&=\sum\limits_{j\in\mathbb{N}}|\langle M_{\sigma}^{*}f,e_j\rangle|^2\\
&=\|M_{\sigma}^{*}f\|^2\\
&=\langle M_{\sigma}M_{\sigma}^{*}f,f\rangle \ \text{for all} \ f\in\mathcal{H}.
\end{align*}
This gives $AKK^{*}\leq M_{\sigma}M_{\sigma}^{*}$.

\vspace{5pt}

$(ii) \Rightarrow (i)$ Let $\{\sigma_i\}_{i\in[m]}$ be any partition of $\mathbb{N}$. Then, by using $(ii)$, for all $f\in\mathcal{H}$ we have
\begin{align*}
&A\|K^{*}f\|^2\\
&=A\langle KK^{*}f,f\rangle\\
&\leq\langle M_{\sigma}M_{\sigma}^{*}f,f\rangle\\
&=\|M_{\sigma}^{*} f\|^2\\
&=\sum\limits_{j\in\mathbb{N}}|\langle M_{\sigma}^{*} f, e_j\rangle|^2\\
&=\sum\limits_{i\in[m]}\sum\limits_{j\in \sigma_{i}}|\langle f,\phi_{ij}\rangle|^2.
\end{align*}
This gives the lower $K$-frame inequality. On the other hand, by Proposition \ref{pro1} the positive number $\sum_{i\in[m]}B_i$ is one of the choice of an universal upper $K$-frame bound. Hence, the family $\Big\{\{\phi_{ij}\}_{j=1}^\infty:i\in[m]\Big\}$   is $K$-woven.
\endproof
Next we give an applicative example of Theorem \ref{th2}.
\begin{exa}
Let $\mathcal{H}=\ell^2(\mathbb{N})$, $m=2$ and $\{e_j\}_{j =1}^{\infty}$ be the canonical orthonormal basis for $\mathcal{H}$.
\begin{enumerate}[(a)]
\item  Define families $\Phi \equiv \{\phi_{1j}\}_{j=1}^\infty$ and $\Psi \equiv \{\phi_{2j}\}_{j=1}^\infty$ in  $\mathcal{H}$ as follows:
\begin{align*}
\{\phi_{1j}\}_{j=1}^\infty&=\{0, e_2, 0, e_3, 0, e_4, 0, \dots \},\\
\{\phi_{2j}\}_{j=1}^\infty&=\{0, e_2, e_2, e_3, e_3, e_4, e_4, \dots \}.
\end{align*}
Let $K$ be the orthogonal projection of $\mathcal{H}$ onto $[e_j]_{j=2}^{\infty}$.
For any $\sigma\subset \mathbb{N}$, define a bounded linear operator $M_{\sigma}:\ell^2(\mathbb{N})\rightarrow\mathcal{H}$ as follows:
\begin{align*}
M_{\sigma}(e_j)= &\begin{cases}
\phi_{1j}, j\in\sigma,\\
\phi_{2j}, j\in\sigma^c.
\end{cases}
\end{align*}
Then, for a given  $f= \{\alpha_1,\alpha_2, \dots\} \in \mathcal{H}$ we have $M^*_{\sigma}(f)=(0,\overline{\alpha_2},\widetilde{\alpha_2},\overline{\alpha_3},\widetilde{\alpha_3},\dots)$, where $\overline{\alpha_n}$ the complex conjugate of the scalar $\alpha_n$ and  $\widetilde{\alpha_n}$ is given by
\begin{align*}
\widetilde{\alpha_n}= &\begin{cases}
0, \ 2n-1\in\sigma \quad (n \in \mathbb{N}),\\
\overline{\alpha_n}, \ 2n-1\in\sigma^c \quad (n \in \mathbb{N}).
\end{cases}
\end{align*}
Thus, for  any $f= \{\alpha_1,\alpha_2, \dots\} \in \mathcal{H}$, we have
\begin{align*}
\langle  M_{\sigma}M_{\sigma}^{*}f , f\rangle=\|M^*_{\sigma}f\|^2\geq\sum\limits_{i\geq2}|\overline{\alpha_i}|^2=\|K^*f\|^2=\langle KK^{*}f, f\rangle.
 \end{align*}
Hence,  by Theorem \ref{th2}, $\Phi$ and  $\Psi$  are $K$-woven.

\vspace{8pt}

\item Define families $\Phi_1 \equiv \{\phi_{1j}\}_{j=1}^\infty$ and $\Psi_1\equiv \{\phi_{2j}\}_{j=1}^\infty$ in  $\mathcal{H}$ as follows:
\begin{align*}
\{\phi_{1j}\}_{j=1}^\infty&=\{e_1, e_2, 0, e_3, e_4, e_5, \dots\},\\
\{\phi_{2j}\}_{j=1}^\infty&=\{e_1, 0, e_2, e_3, e_4, e_5, \dots\}.
\end{align*}
Let $K$ be the orthogonal projection of $\mathcal{H}$ onto $[e_j]_{j=2}^{\infty}$. Then, $\Phi_1$ and  $\Psi_1$ are $K$-frames for $\mathcal{H}$.

 Next we show that $\Phi_1$ and $\Psi_1$ are not $K$-woven. Choose  $\sigma=\mathbb{N}\setminus\{2\}$. Assume that   there exists  a bounded linear operator $M_{\sigma}:\ell^2(\mathbb{N})\rightarrow\mathcal{H}$ such that
\begin{align*}
M_{\sigma}(e_j)= &\begin{cases}
\phi_{1j}, \  j\in\sigma,\\
\phi_{2j}, \ j\in\sigma^c.
\end{cases}
\end{align*}
Then, for any $A>0$, we have
\begin{align*}
A\langle KK^{*}e_2, e_2\rangle &=A\langle K^{*}e_2, K^{*} e_2\rangle\\
 &=A\|e_2\|^2\\
& > 0 \\
&= \|M_{\sigma} e_2\|^2\\
& = \langle  M_{\sigma}M_{\sigma}^{*}e_2 ,e_2\rangle.
\end{align*}
Hence,  by Theorem  \ref{th2},  $\Phi_1$ and $\Psi_1$ are  not $K$-woven.
\end{enumerate}
\end{exa}

 The following theorem provides a  sufficient condition for  $K$-frames  not to be weakly woven  in terms of a lower $K$-frame bound. This is inspired by   \cite[Lemma 4.3]{BCGLL}.

\begin{thm}\label{th2.7}
Suppose $\{\phi_{ij}\}_{j=1}^\infty$ is a K-frame for $\mathcal{H}$ with bounds $A_i$ and $B_i$  $(i\in [m])$. Assume for any partition  $\{\tau_i\}_{i\in [m]}$ of a finite subset of $\mathbb{N}$ and for every $A>0$ there exists  a partition $\{\sigma_i\}_{i\in[m]}$ of $\mathbb{N}\setminus\{\tau_i\}_{i\in [m]}$ such that $\bigcup_{i\in[m]}\{\phi_{ij}\}_{j\in\sigma_i\cup\tau_i}$ has a lower $K$-frame bound less than $A$. Then, there exists  a partition $\{\pi_i\} _{i\in[m]}$ of $\mathbb{N}$ such that $\bigcup_{i\in[m]}\{\phi_{ij}\}_{j\in\pi_i}$ is not a $K$-frame for $\mathcal{H}$.
\end{thm}
\proof
Let $\tau_{1i}= \varnothing$ for all $i\in[m]$. Then,  for $A_{1}=1$, there exists a partition $\{\sigma_{1i}\}_{i\in[m]}$ of $\mathbb{N}$ such that $\bigcup_{i\in[m]}\{\phi_{ij}\}_{j\in\sigma_{1i}\cup\tau_{1i}}$ has a lower $K$-frame  bound less than $1$. Therefore, there exists a vector $h_1\in \mathcal{H}$ such that
\begin{align*}
\sum\limits_{i\in[m]}&\sum\limits_{j\in \sigma_{1i}\cup\tau_{1i}}|\langle h_1,\phi_{ij}\rangle|^2<\|K^{*}h_1\|^2.
\end{align*}
Since $\sum\limits_{i\in[m]}\sum\limits_{j\in \mathbb{N}}|\langle h_1,\phi_{ij}\rangle|^2<\infty$,  there exists  $k_1\in \mathbb{N}$ such that
\begin{align*}
\sum\limits_{i\in[m]}\sum\limits_{j>k_1}|\langle h_1,\phi_{ij}\rangle|^2<\|K^{*}h_1\|^2.
\end{align*}
Choose a partition $\{\tau_{2i}\}_{i\in[m]}$ of $[k_1]$ such that $\tau_{2i}=\tau_{1i}\cup(\sigma_{1i}\cap[k_1])$ for all $i\in[m]$ and $A_{2}=\frac{1}{2}$. Then,  there exists  a partition $\{\sigma_{2i}\}_{i\in[m]}$ of $\mathbb{N}\setminus\{\tau_{2i}\}_{i\in [m]}$ such that the family  $\bigcup_{i\in[m]}\{\phi_{ij}\}_{j\in\sigma_{2i}\cup\tau_{2i}}$ has a lower $K$-frame  bound less than $\frac{1}{2}$. Therefore, there exists a vector $h_2\in \mathcal{H}$ such that
\begin{align*}
\sum\limits_{i\in[m]}\sum\limits_{j\in \sigma_{2i}\cup\tau_{2i}}|\langle h_2,\phi_{ij}\rangle|^2<\frac{1}{2}\|K^{*}h_2\|^2.
\end{align*}
By hypothesis, we have
\begin{align*}
\sum\limits_{i\in[m]}\sum\limits_{j\in \mathbb{N}}|\langle h_2,\phi_{ij}\rangle|^2<\infty.
\end{align*}
Therefore,  there exists  $k_2>k_1$ such that
\begin{align*}
\sum\limits_{i\in[m]}\sum\limits_{j>k_2}|\langle h_2,\phi_{ij}\rangle|^2<\frac{1}{2}\|K^{*} h_2\|^2.
\end{align*}
Continuing in this way, for $A_{n}=\frac{1}{n}$ and for a partition $\{\tau_{ni}\}_{i\in[m]}$ of $[k_{n-1}]$ such that \break  $\tau_{ni}=\tau_{(n-1)i}\cup(\sigma_{(n-1)i}\cap[k_{n-1}])$  $(i\in[m])$, there exists  a partition $\{\sigma_{ni}\}_{i\in[m]}$ of $\mathbb{N}\setminus\{\tau_{ni}\}_{i\in [m]}$ such that $\bigcup_{i\in[m]}\{\phi_{ij}\}_{j\in\sigma_{ni}\cup\tau_{ni}}$ has a lower $K$-frame bound less than $\frac{1}{n}$. Thus,  there exists  $h_n\in \mathcal{H}$ such that
\begin{align}\label{eq3.3}
\sum\limits_{i\in[m]}\sum\limits_{j\in \sigma_{ni}\cup\tau_{ni}}|\langle h_n,\phi_{ij}\rangle|^2<\frac{1}{n}\|K^{*}h_n\|^2,
\end{align}
and we can find a positive integer  $k_n>k_{n-1}$ such that
\begin{align}\label{eq3.4}
\sum\limits_{i\in[m]}\sum\limits_{j>k_n}|\langle h_n,\phi_{ij}\rangle|^2<\frac{1}{n}\|K^{*}h_n\|^2.
\end{align}

Choose a partition $\{\pi_i\} _{i\in[m]}$ of $\mathbb{N}$, where  $\pi_i=\bigcup_{n\in\mathbb{N}}\{\tau_{ni}\}$. We show that the  family
 $\bigcup_{i\in[m]}\{\phi_{ij}\}_{j\in\pi_i}$ is not a $K$-frame for $\mathcal{H}$. Assume that   $\bigcup_{i\in[m]}\{\phi_{ij}\}_{j\in\pi_i}$ is a  $K$-frame for $\mathcal{H}$ with bounds $\alpha$ and $\beta$.  By  the  Archimedean property there exists a  $\eta\in \mathbb{N}$ such that $\eta>\frac{2}{\alpha}$.

  Using \eqref{eq3.3} and  \eqref{eq3.4},  we compute
\begin{align*}
&\sum_{i\in[m]}\sum_{j\in \pi_i}|\langle h_\eta,\phi_{ij}\rangle|^2\\
&\leq \sum_{i\in[m]} \sum_{j\in \tau_{(\eta+1)i}}|\langle h_\eta,\phi_{ij}\rangle|^2+\sum_{i\in[m]} \sum_{j\geq k_\eta}|\langle h_\eta,\phi_{ij}\rangle|^2\\
&\leq \sum_{i\in[m]} \sum_{j\in \tau_{\eta i}\cup\sigma_{\eta i}}|\langle h_\eta,\phi_{ij}\rangle|^2+\sum_{i\in[m]} \sum_{j\geq k_\eta}|\langle h_\eta,\phi_{ij}\rangle|^2\\
&\leq \frac{1}{\eta}\|K^{*}h_\eta\|^2+\frac{1}{\eta}\|K^{*} h_\eta\|^2\\
&=\frac{2}{\eta}\|K^{*}h_\eta\|^2\\
&<\alpha\|K^{*}h_\eta\|^2, \ \text{a contradiction }.
\end{align*}
This completes the proof.
\endproof
Theorem \ref{th2.7} gives a necessary condition for weakly woven $K$-frames.
\begin{prop}\label{pro2.9}
Suppose the family of $K$-frames $\Big\{\{\phi_{ij}\}_{j=1}^\infty:i\in[m]\Big\}$ for $\mathcal{H}$ is weakly $K$-woven. Then, there exist a partition  $\{\tau_i\}_{i\in[m]}$ of some finite subset of $\mathbb{N}$ and $A>0$ such that for any partition $\{\sigma_i\}_{i\in[m]}$ of $\mathbb{N}\setminus\{\tau_i\}_{i\in[m]}$,   the family $\bigcup\limits_{i\in[m]}\{\phi_{ij}\}_{j\in\sigma_i\cup\tau_i}$ has a lower $K$-frame bound $A$.
\end{prop}

Next we show that woven $K$-frames are equivalent to weakly  woven  $K$-frames.
\begin{thm}\label{th2.6n}
Suppose $\Phi \equiv \{\phi_{j}\}_{j=1}^\infty$ and $\Psi \equiv \{\psi_{j}\}_{j=1}^\infty$ are $K$-frames for $\mathcal{H}$. The following are equivalent:
\begin{enumerate}[$(i)$]
\item $\Phi$ and $\Psi$ are $K$-woven.
\item $\Phi$ and $\Psi$  are weakly $K$-woven.

\end{enumerate}
\end{thm}
\proof $(i) \Rightarrow (ii):$ Obvious.

$(ii) \Rightarrow (i):$
First we note that an  universal upper $K$-frame bound for $\Phi$ and $\Psi$ can be obtained from Proposition \ref{pro1}. So it is sufficient to compute an  universal lower $K$-frame bound for $\Phi$ and $\Psi$.

By Proposition \ref{pro2.9}, there exist disjoint finite sets $I$ and $J$ of $\mathbb{N}$ and $A>0$ satisfying:

\vspace{5pt}
$(\maltese)$ \quad For any partition $\{\sigma, \delta\}$ of $\mathbb{N}\setminus(I\cup J)$, the family   $\{\phi_{j}\}_{j\in I\cup\sigma}\bigcup\{\psi_{j}\}_{j\in J\cup\delta}$ has a lower $K$-frame bound $A$.

\vspace{5pt}
We can permute both the $K$-frames $\Phi$ and $\Psi$ (if necessary), so that $I\cup J=[m]$.

\vspace{5pt}
\textbf{Claim :} For any  partition $\{I_o, J_o\}$ of $[m]$, there exists  $A_o >0$ such that for any partition $\{\sigma, \delta\}$ of $\mathbb{N}\setminus[m]$ the family $\{\phi_{j}\}_{j\in I_o\cup\sigma}\bigcup\{\psi_{j}\}_{j\in J_o\cup\delta}$ has a lower $K$-frame bound $A_o$. So that $\Phi$ and $\Psi$ are $K$-woven  with an universal lower $K$-frame bound
$\min\Big\{{A_o:\{I_0,J_0\} \text{ is a partition of} \ [m]}\Big\}>0$.

Suppose our claim is not true. Then,  there exists a partition $\{I_1, J_1\}$ of $[m]$ such that for each positive $\alpha$, there exists a partition $\{\sigma_\alpha, \delta_\alpha\}$ of $\mathbb{N}\setminus[m]$ such that $\{\phi_{j}\}_{j\in I_1\cup\sigma_\alpha}\bigcup\{\psi_{j}\}_{j\in J_1\cup\delta_\alpha}$ has a lower $K$-frame bound less than $\alpha$.
Therefore, for any $n\in\mathbb{N}$, there exists a partition $\{\sigma_n, \delta_n\}$ of $\mathbb{N}\setminus[m]$ such that $\{\phi_{j}\}_{j\in I_1\cup\sigma_n}\bigcup\{\psi_{j}\}_{j\in J_1\cup\delta_n}$ has a lower $K$-frame bound less than $n$. Thus, there exists $h_n\in\mathcal{H}$ such that
\begin{align*}
&\sum\limits_{j\in I_1\cup\sigma_{n}}|\langle h_n,\phi_{j}\rangle|^2+\sum\limits_{j\in J_1\cup\delta_{n}}|\langle h_n,\psi_{j}\rangle|^2\\
 & <\frac{1}{n}\|K^{*}(h_n)\|^2\\
&\leq \frac{1}{n}\lambda^2\|h_n\|^2, \text{ where } \lambda^2=\|K^*\|^2.
\end{align*}
That is
\begin{align*}
\sum\limits_{j\in I_1\cup\sigma_{n}}|\langle \zeta g_n,\phi_{j}\rangle|^2+\sum\limits_{j\in J_1\cup\delta_{n}}|\langle \zeta g_n,\psi_{j}\rangle|^2<\frac{1}{n},
\end{align*}
 where $\zeta = \frac{1}{\lambda} \ \text{and} \ g_n = \frac{h_n}{\|h_n\|} \ (n \in \mathbb{N})$.
Note  that $K^*(\zeta g_n)\neq 0 \ (n \in \mathbb{N})$.

\vspace{5pt}
Now we proceed in the following steps:
\vspace{5pt}

\textbf{Step 1:}  By hypothesis, for each $n\in \mathbb{N}$ there exist a partition $\{\sigma_n, \delta_n\}$ of $\mathbb{N}\setminus[m]$ and a unit vector $g_n\in\mathcal{H}$ such that
\begin{align}\label{eq2.4}
\sum\limits_{j\in I_1\cup\sigma_{n}}|\langle \zeta g_n,\phi_{j}\rangle|^2+\sum\limits_{j\in J_1\cup\delta_{n}}|\langle \zeta g_n,\psi_{j}\rangle|^2<\frac{1}{n}
\end{align}
and the sets $\sigma_n$ and $\delta_n$ satisfy the following properties:
\begin{enumerate}
\item For every $k=1,2, \dots$, either $m+k\in\sigma_n$ for all $n\geq k$ or $m+k\in\delta_n$ for all $n\geq k$.\label{eqint1}
\item There is a partition $\{\sigma, \delta\}$ of $\mathbb{N}\setminus[m]$ such that $m+k\in\sigma$ implies that $m+k\in\sigma_n$ for all $n\geq k$ or if $m+k\in\delta$ implies that $m+k\in\delta_n$ for all $n\geq k$.\label{eqint2}
\end{enumerate}

Furthermore, the sequence  $\{\zeta g_n\}_{n=1}^{\infty}$ is bounded, so there exists a subsequence $\{\zeta g_{n_i}\}_{i=1}^{\infty}$ of $\{\zeta g_n\}_{n=1}^{\infty}$ such that  $\{\zeta g_{n_i}\}_{i=1}^{\infty}$ converges weakly to $h \in \mathcal{H}$. We reindex, $\zeta g_{n_i}\rightarrow\zeta g_i$ and $\sigma_{n_i}\rightarrow\sigma_i$, $\delta_{n_i}\rightarrow\delta_i$. Note that  \eqref{eq2.4},  \eqref{eqint1} and \eqref{eqint2} are  satisfied by this constructed sequence.

\vspace{10pt}
\textbf{Step 2:} In this step we show that $K^*(h)\neq 0$. Let $\widetilde{K^{*}}:N(K^{*})^{\bot}\rightarrow R(K^*)$ be the restriction of $K^{*}$ on $N(K^{*})^{\bot}$. Note that $R(K^*)$ is closed, $\widetilde{K^{*}}$ is invertible and  $\widetilde{K^*}^{-1}$ is bounded.\\
For any non-zero  $f\in N(K^{*})^{\bot}$,  $f=\widetilde{K^*}^{-1}(g)$ for some $g\in R(K^*)$.

We compute
\begin{align}\label{eq2.5}
\frac{\|K^* f\|}{\|f\|}&=\frac{\|\widetilde{K^*} f\|}{\|f\|}
=\frac{\|\widetilde{K^*}(\widetilde{K^*}^{-1} g)\|}{\|\widetilde{K^*}^{-1} g \|}
=\frac{\|g\|}{\|\widetilde{K^*}^{-1} g\|}
\geq\frac{1}{\|\widetilde{K^*}^{-1}\|}.
\end{align}
Fix   $p\in\mathbb{N}$ such that $p>\frac{2\lambda^2\|\widetilde{K^*}^{-1}\|^2}{A}$. So, we can find $n_p\in\mathbb{N}$ such that for all $n\geq n_p>p$, we have
\begin{align}\label{eq2.6}
\sum\limits_{j\in [m+p]}|\langle \zeta g_n-h,\phi_{j}\rangle|^2+\sum\limits_{j\in[m+p]}|\langle \zeta g_n-h,\psi_{j}\rangle|^2<\frac{1}{2p}.
\end{align}
From  (\maltese), for each $n\in\mathbb{N}$ we have
\begin{align}\label{eq2.7}
\sum\limits_{j\in I\cup\sigma_{n}}|\langle \zeta g_n,\phi_{j}\rangle|^2+\sum\limits_{j\in J\cup\delta_{n}}|\langle \zeta g_n,\psi_{j}\rangle|^2\geq A\|K^*(\zeta g_n)\|^2.
\end{align}
For $m, p \in \mathbb{N}$, we write $[m, p] = \{m+1, m+2, \dots, m+p\}$. By using \eqref{eq2.4}, \eqref{eq2.5}, \eqref{eq2.6}, \eqref{eq2.7} and $\sigma\cap[m,p]=\sigma_{n_p}\cap[m,p]$, $\delta\cap[m,p]=\delta_{n_p}\cap[m,p]$, we compute
\begin{align*}
&\sum\limits_{j\in I\cup\sigma}|\langle h,\phi_{j}\rangle|^2+\sum\limits_{j\in J\cup\delta}|\langle h,\psi_{j}\rangle|^2\\
&\geq\sum\limits_{j\in I\cup(\sigma\cap[m,p])}|\langle h,\phi_{j}\rangle|^2+\sum\limits_{j\in J\cup(\delta\cap[m,p])}|\langle h,\psi_{j}\rangle|^2\\
&=\sum\limits_{j\in I\cup(\sigma_{n_p}\cap[m,p])}|\langle h,\phi_{j}\rangle|^2+\sum\limits_{j\in J\cup(\delta_{n_p}\cap[m,p])}|\langle h,\psi_{j}\rangle|^2\\
&\geq\frac{1}{2}\left(\sum\limits_{j\in I\cup(\sigma_{n_p}\cap[m,p])}|\langle \zeta g_{n_p},\phi_{j}\rangle|^2+\sum\limits_{j\in J\cup(\delta_{n_p}\cap[m,p])}|\langle \zeta g_{n_p},\psi_{j}\rangle|^2\right)\\
& -\left(\sum\limits_{j\in I\cup(\sigma_{n_p}\cap[m,p])}|\langle h-\zeta g_{n_p},\phi_{j}\rangle|^2+\sum\limits_{j\in J\cup(\delta_{n_p}\cap[m,p])}|\langle h-\zeta g_{n_p},\psi_{j}\rangle|^2\right)\\
&\geq\frac{1}{2}\left(\sum\limits_{j\in I\cup\sigma_{n_p}}|\langle \zeta g_{n_p},\phi_{j}\rangle|^2+\sum\limits_{j\in J\cup\delta_{n_p}}|\langle \zeta g_{n_p},\psi_{j}\rangle|^2\right)\\
& -\frac{1}{2}\left(\sum\limits_{j\in \sigma_{n_p}\cap[m+p]^c}|\langle \zeta g_{n_p},\phi_{j}\rangle|^2+\sum\limits_{j\in \delta_{n_p}\cap[m+p]^c}|\langle \zeta g_{n_p},\psi_{j}\rangle|^2\right)\\
& -\left(\sum\limits_{j\in [m+p]}|\langle h-\zeta g_{n_p},\phi_{j}\rangle|^2+\sum\limits_{j\in [m+p]}|\langle h-\zeta g_{n_p},\psi_{j}\rangle|^2\right)\\
&\geq \frac{A}{2}\|K^*(\zeta g_{n_p})\|^2-\frac{1}{2}\frac{1}{n_p}-\frac{1}{2p}\\
&= \frac{A\zeta^2}{2}\|K^*( g_{n_p})\|^2-\frac{1}{2}\frac{1}{n_p}-\frac{1}{2p}\\
&\geq \frac{A}{2\lambda^2\|\widetilde{K^*}^{-1}\|^2}-\frac{1}{p}\\
&>0.
\end{align*}
This gives $h\neq 0$. Since $\{\zeta g_n\}_{n = 1}^{\infty}\subset  N(K^{*})^{\bot}$ and $N(K^{*})^{\bot}$ is  closed, so  $K^*(h)\neq 0$.

\vspace{10pt}
\textbf{Step 3:} In this step we will show that $\{\phi_j\}_{j\in I_1\cup\sigma} \bigcup \{\psi_j\}_{j\in J_1\cup\delta}$ is not a $K$-frame for $\mathcal{H}$.\\
By using \eqref{eq2.4} and \eqref{eq2.6}, we compute
\begin{align*}
&\sum\limits_{j\in I_1\cup\sigma}|\langle h,\phi_{j}\rangle|^2+\sum\limits_{j\in J_1\cup\delta}|\langle h,\psi_{j}\rangle|^2\\
&=\lim_{p\rightarrow\infty}\left(\sum\limits_{j\in I_1\cup(\sigma\cap[m,p])}|\langle h,\phi_{j}\rangle|^2+\sum\limits_{j\in J_1\cup(\delta\cap[m,p])}|\langle h,\psi_{j}\rangle|^2\right)\\
&=\lim_{p\rightarrow\infty}\left(\sum\limits_{j\in I_1\cup(\sigma_{n_p}\cap[m,p])}|\langle h,\phi_{j}\rangle|^2+\sum\limits_{j\in J_1\cup(\delta_{n_p}\cap[m,p])}|\langle h,\psi_{j}\rangle|^2\right)\\
&\leq2\lim_{p\rightarrow\infty}\left(\sum\limits_{j\in I_1\cup(\sigma_{n_p}\cap[m,p])}|\langle \zeta g_{n_p},\phi_{j}\rangle|^2+\sum\limits_{j\in J_1\cup(\delta_{n_p}\cap[m,p])}|\langle \zeta g_{n_p},\psi_{j}\rangle|^2\right)\\
& +2\lim_{p\rightarrow\infty}\left(\sum\limits_{j\in I_1\cup(\sigma_{n_p}\cap[m,p])}|\langle h-\zeta g_{n_p},\phi_{j}\rangle|^2+\sum\limits_{j\in J_1\cup(\delta_{n_p}\cap[m,p])}|\langle h-\zeta g_{n_p},\psi_{j}\rangle|^2\right)\\
&\leq2\lim_{p\rightarrow\infty}\frac{1}{n_p}+2\lim_{p\rightarrow\infty}\frac{1}{2p}\\
&=0\\
& <\alpha\|K^* h\|^2 \text{ for any } \alpha>0.
\end{align*}
Therefore, $\{\phi_{j}\}_{j\in I_1\cup\sigma}\bigcup\{\psi_{j}\}_{j\in J_1\cup\delta}$ is not a $K$-frame for $\mathcal{H}$. This  contradicts the fact that  $\Phi$ and $\Psi$ are weakly $K$-woven. The proof is complete.
\endproof

\begin{rem}\label{rem9ext1}
The technique used in the proof of Theorem  \ref{th2.6n} is different than  given in  \cite[Theorem 4.5]{BCGLL} for an analogous result for ordinary frames in separable Hilbert spaces. In the proof of Theorem  \ref{th2.6n},  the operator $K$ is not invertible on the whole space $\mathcal{H}$. So, computations in our proof are different than given in
\cite[Theorem 4.5]{BCGLL}. Furthermore,    if $K = I$, the identity operator on $\mathcal{H}$, then by Theorem \ref{th2.6n}, we
can  obtain \cite[Theorem 4.5]{BCGLL}.
\end{rem}

The following result  characterizes weaving $K$-frames in terms of action of operators on  $K$-frames.
\begin{prop}\label{pr2}
Let  $\Big\{\{\phi_{ij}\}_{j=1}^\infty:i\in[m]\Big\}$ be a family of $K$-frames for $\mathcal{H}$. The following are equivalent:
\begin{enumerate}[$(i)$]
\item  $\Big\{\{\phi_{ij}\}_{j=1}^\infty:i\in[m]\Big\}$ is $K$-woven.
\item $\Big\{\{U(\phi_{ij})\}_{j=1}^\infty:i\in[m]\Big\}$ is $UK$-woven for all $U \in \mathcal{B}(\mathcal{H})$.
\end{enumerate}
\end{prop}
\proof
$(i) \Rightarrow (ii):$ Let $A$ and $B$ be universal $K$-frame bounds for the family $\Big\{\{\phi_{ij}\}_{j=1}^\infty:i\in[m]\Big\}$. Let
$\{\sigma_{i}\}_{i\in [m]}$ be  any partition of $\mathbb{N}$. Then, for any $f\in\mathcal{H}$, we have
\begin{align*}
\sum\limits_{i\in[m]}\sum\limits_{j\in \sigma_{i}}|\langle f,U(\phi_{ij})\rangle|^2=\sum\limits_{i\in[m]}\sum\limits_{j\in \sigma_{i}}|\langle U^{*}f,\phi_{ij}\rangle|^2
\leq B\|U^{*}f\|^2
\leq B\|U^{*}\|^2\|f\|^2.
\end{align*}
Similarly, for any $f \in \mathcal{H}$, we have
\begin{align*}
\sum\limits_{i\in[m]}\sum\limits_{j\in \sigma_{i}}|\langle f,U(\phi_{ij})\rangle|^2=\sum\limits_{i\in[m]}\sum\limits_{j\in \sigma_{i}}|\langle U^{*} f,\phi_{ij}\rangle|^2
\geq A\|K^{*}U^{*} f\|^2
= A\|(UK)^{*}f\|^2.
\end{align*}
Hence, the family $\Big\{\{U(\phi_{ij})\}_{j=1}^\infty:i\in[m]\Big\}$ is $UK$-woven with universal $K$-frame bounds $A, \ B\|U^{*}\|^2$.

$(ii) \Rightarrow (i):$ Choose $U = I$, the identity operator on $\mathcal{H}$. Then, the family $\Big\{\{\phi_{ij}\}_{j=1}^\infty:i\in[m]\Big\}$ is $K$-woven.
\endproof
\begin{rem}
Let  $\Phi$ and $\Psi$ be  $K$-frames for $\mathcal{H}$ such that  $U \Phi$ and $U \Psi$ are $UK$-woven for some  $U \in \mathcal{B}(\mathcal{H})$. Then, in general, $\Phi$ and $\Psi$ are not $K$-woven. This is justified in the following example.
\end{rem}
\begin{exa}
Let $\mathcal{H}=\ell^2(\mathbb{N})$ and let $K$ be the orthogonal projection of $\mathcal{H}$ onto $[e_j]_{j=2}^{\infty}$.
Define $\Phi, \Psi \subset \mathcal{H}$ by
\begin{align*}
\Phi \equiv \{\phi_{1j}\}_{j=1}^\infty&=\{e_1, 0, e_2, 0, e_3, 0, e_4, 0, e_5, 0 \dots \},\\
\Psi \equiv \{\phi_{2j}\}_{j=1}^\infty&=\{0, e_1, 0, e_2, e_3, e_3, e_4, e_4, e_5, e_5, \dots\},
\end{align*}
where $\{e_j\}_{j=1}^{\infty}$ is the canonical orthonormal basis for $\mathcal{H}$. Then, $\Phi$ and $\Psi$ are  $K$-frames for $\mathcal{H}$.

 To show $\Phi$ and $\Psi$ are not $K$-woven. Choose $\sigma=\mathbb{N}\setminus\{1,3\}$. Then, the family  $\{\phi_{1j}\}_{j\in\sigma}\bigcup\{\phi_{2j}\}_{j\in\sigma^c}=\{0,0,0,0,e_3,0,e_4,0,e_5,0,e_6, \dots\}$ is not a $K$-frame  for $\mathcal{H}$, since for any $A>0$, we have
\begin{align*}
&\sum_{j\in\sigma}|\langle e_2,\phi_{1j}\rangle|^2+\sum_{j\in\sigma^c}|\langle e_2,\phi_{2j}\rangle|^2\\
 &=\sum_{j\geq 3}|\langle e_2,e_j\rangle|^2\\
&=0\\
&< A\|e_2\|^2\\
&=A\|K^*e_2\|^2.
\end{align*}
 Hence $\Phi$ and  $\Psi$ are not $K$-woven.

  Let $U$ be the orthogonal projection of $\mathcal{H}$ onto $[e_j]_{j=3}^{\infty}$.
Then,   $U \Phi \equiv \{U(\phi_{1j})\}_{j=1}^\infty$ and  $U \Psi \equiv  \{U(\phi_{2j})\}_{j=1}^\infty$ are $UK$-woven frames for $\mathcal{H}$.  To show this, first we note that:
\begin{align*}
\{U(\phi_{1j})\}_{j=1}^\infty&=\{0, 0, 0, 0, e_3, 0, e_4, 0, e_5, 0, e_6, \dots\},\\
\{U(\phi_{2j})\}_{j=1}^\infty&=\{0,0,0,0,e_3,e_3,e_4,e_4,e_5,e_5,e_6, \dots\}.
\end{align*}
For any subset $\sigma$ of $\mathbb{N}$, we have
\begin{align*}
\sum_{j\in\sigma}|\langle f, U(\phi_{1j})\rangle|^2+\sum_{j\in\sigma^c}|\langle f, U(\phi_{2j})\rangle|^2\leq2\sum_{j=1}^{\infty}|\langle f,e_j\rangle|^2=2\|f\|^2, \ f\in\mathcal{H}.
\end{align*}
On the other hand, let $f\in\mathcal{H}$. Then,  $f=\sum\limits_{j = 1}^{\infty}\alpha_je_j$. Thus, we have
\begin{align*}
\|(UK)^{*}f\|^2
&=\|U^{*}f\|^2\\
&=\|\sum\limits_{j\in\mathbb{N}}\alpha_jU^{*}e_j\|^2\\
&=\|\sum\limits_{j=3}^{\infty}\alpha_je_j\|^2\\
&=\sum\limits_{j=3}^{\infty}|\langle f,e_j\rangle|^2\\
&\leq \sum_{j\in\sigma}|\langle f,U\phi_{1j}\rangle|^2+\sum_{j\in\sigma^c}|\langle f, U\phi_{2j}\rangle|^2.
\end{align*}
Hence  $U \Phi$ and $U \Psi$ are $UK$-woven with universal $UK$-frame bounds $1,  2$.
\end{exa}

The next result gives a  sufficient conditions for perturbation of weaving $K$-frames.
\begin{thm}\label{th8}
Suppose $\Phi_1 \equiv\{\phi_{1j}\}_{j=1}^\infty$  and $\Phi_2 \equiv\{\phi_{2j}\}_{j=1}^\infty$  are $K$-frames  for $\mathcal{H}$ with bounds $A_1, \ B_1$ $($respectively, $A_2, \ B_2$$)$ such that $\Phi_1$ is an orthogonal set and  $\|\phi_{1j}\|^2>\alpha>0$ for all $j\in\mathbb{N}$. Assume that there exists $\lambda, \mu, \nu\geq0$ such that
\begin{align}\label{eq2.7asha}
\Big(\sqrt{B_1}+\sqrt{B_2}\Big)\left(\lambda+\mu\sqrt{B_1}+\nu\sqrt{B_2}\right)\|\widetilde{K^{*}}^{-1}\| < \sqrt{\alpha A_1}
\end{align}
and
\begin{align}\label{2.5}
\Big\|\sum_{j=1}^{\infty} a_j(\phi_{1j}-\phi_{2j})\Big\|\leq \lambda \|\{a_j\}_{j=1}^{\infty}\|_{\ell^2(\mathbb{N})}+\mu\Big\|\sum_{j=1}^{\infty} a_j\phi_{1j}\Big\|+\nu\Big\|\sum_{j=1}^{\infty}a_j\phi_{2j}\Big\|
\end{align}
for all sequence of scalars $\{a_j\}_{j=1}^{\infty} \in  \ell^2(\mathbb{N})$, where $\widetilde{K^{*}}:N(K^{*})^{\bot}\rightarrow R(K^*)$ is the restriction of $K^{*}$ on $N(K^{*})^{\bot}$. Then,  $\Phi_1$  and $\Phi_2$ are $K$-woven with universal  bounds \break $\frac{\Big[\sqrt{\alpha A_1} - (\sqrt{B_1}+\sqrt{B_2})(\lambda+\mu\sqrt{B_1}+\nu\sqrt{B_2})\|\widetilde{K^{*}}^{-1}\|\Big]^2}{B_1 + B_2} $ and $B_1 + B_2$.
\end{thm}
\proof

Let  $T^i$ be the pre-frame operator associated with $\Phi_i, \ i = 1, 2$. Then, $\|T^i \| \leq \sqrt{B_i}$.
Let $\sigma$ be any subset of $\mathbb{N}$. Define  $T^i_{\sigma^c}: \ell^2(\mathbb{N}) \rightarrow \mathcal{H}$ by
\begin{align*}
T^i_{\sigma^c}(\{a_j\}_{j=1}^{\infty})=\sum\limits_{j\in \sigma^c}a_j\phi_{ij}, \ \{a_j\}_{j=1}^{\infty} \in \ell^2(\mathbb{N}) \ \quad (i =1, 2).
\end{align*}
Then, each $T^i_{\sigma^c}$  is a bounded operator with norm at most  $\sqrt{B_i}$. Indeed, for all $\{a_j\}_{j=1}^{\infty} \in \ell^2(\mathbb{N})$, we have
\begin{align*}
\left\|T^i_{\sigma^c}(\{a_j\}_{j\in \mathbb{N}})\right\|
&=\Big\|\sum_{j\in \sigma^c}a_j\phi_{ij}\Big\|\\
&\leq \|T^i\| \|\{a_j\}_{j=1}^{\infty}\|_{\ell^2(\mathbb{N})}\\
&\leq\sqrt{B_i}\|\{a_j\}_{j=1}^{\infty}\|_{\ell^2(\mathbb{N})} \ \quad (i =1, 2).
\end{align*}

By using \eqref{2.5}, we compute
\begin{align}\label{eq2.3}
&\|(T^1_{\sigma^c}-T^2_{\sigma^c})\{a_j\}_{j\in \mathbb{N}}\| \notag\\
&=\Big\|\sum_{j\in\sigma^c}a_j(\phi_{1j}-\phi_{2j})\Big\| \notag\\
&=\Big\|\sum_{j\in \mathbb{N}}\gamma_j(\phi_{1j}-\phi_{2j})\Big\| \ (\text{where} \ \gamma_j=a_j \text{ for } j\in \sigma^c \ \text{and} \  \gamma_j =0 \ \text{otherwise}) \notag\\
&\leq\lambda\left\|\{\gamma_j\}_{j=1}^{\infty}\right\|_{\ell^2(\mathbb{N})}+\mu\Big\|\sum\limits_{j=1}^{\infty}\gamma_j\phi_{1j}\Big\|+\nu\Big\|\sum\limits_{j=1}^{\infty} \gamma_j\phi_{2j}\Big\| \notag\\
&\leq\lambda\left\|\{\gamma_j\}_{j=1}^{\infty}\right\|_{\ell^2(\mathbb{N})}+\mu\|T^1\|\left\|\{\gamma_j\}_{j=1}^{\infty}\right\|_{\ell^2(\mathbb{N})}+\nu\|T^2\|\left\|\{\gamma_j\}_{j=1}^{\infty}\right\|_{\ell^2(\mathbb{N})}\notag\\
&\leq\left(\lambda+\mu\sqrt{B_1}+\nu\sqrt{B_2}\right)\|\{a_j\}_{j=1}^{\infty}\|_{\ell^2(\mathbb{N})} \  \text{ for all } \{a_j\}_{j=1}^{\infty} \in \ell^2(\mathbb{N}).
\end{align}

For any $f\in N(K^{*})^{\bot}$, by using \eqref{eq2.3} and orthogonality of $\Phi_1$,  we compute
\begin{align}\label{eq2.11}
&\Big\|\sum_{j\in \sigma}\langle f, \phi_{1j}\rangle\phi_{1j}+\sum_{j\in \sigma^c}\langle f, \phi_{2j}\rangle\phi_{2j}\Big\| \notag\\
&=\Big\|\sum_{j\in \sigma}\langle f, \phi_{1j}\rangle\phi_{1j}+\sum_{j\in \sigma^c}\langle f, \phi_{1j}\rangle\phi_{1j}-\sum_{j\in \sigma^c}\langle f, \phi_{1j}\rangle\phi_{1j}+\sum_{j\in \sigma^c}\langle f, \phi_{2j}\rangle\phi_{2j}\Big\| \notag\\
&=\Big\|\sum_{j\in\mathbb{N}}\langle f, \phi_{1j}\rangle\phi_{1j}-\sum_{j\in \sigma^c}\langle f, \phi_{1j}\rangle\phi_{1j}+\sum_{j\in \sigma^c}\langle f, \phi_{2j}\rangle\phi_{2j}\Big\| \notag\\
&\geq\Big\|\sum_{j\in\mathbb{N}}\langle f, \phi_{1j}\rangle\phi_{1j}\Big\|-\Big\|\sum_{j\in \sigma^c}\langle f, \phi_{1j}\rangle\phi_{1j}-\sum_{j\in \sigma^c}\langle f, \phi_{2j}\rangle\phi_{2j}\Big\| \notag\\
&\geq{\Big(\sum_{j\in\mathbb{N}}|\langle f,\phi_{1j}\rangle|^2\|\phi_{1j}\|^2\Big)}^{\frac{1}{2}}- \Big\|\sum_{j\in \sigma^c}\langle f, \phi_{1j}\rangle\phi_{1j}-\sum_{j\in \sigma^c}\langle f, \phi_{1j}\rangle\phi_{2j}\Big\| \notag\\
& \quad -\Big\|\sum_{j\in \sigma^c}\langle f, \phi_{1j}\rangle\phi_{2j}-\sum_{j\in \sigma^c}\langle f, \phi_{2j}\rangle\phi_{2j}\Big\| \notag\\
&\geq\sqrt{\alpha A_1}\|K^{*}(f)\|-\|T^1_{\sigma^c}T^{1*}(f)-T^2_{\sigma^c}T^{1*}(f)\|-\|T^2_{\sigma^c}T^{1*}(f)-T^2_{\sigma^c}T^{2*}(f)\| \notag\\
&=\sqrt{\alpha A_1}\|K^{*}(f)\|-\|T^1_{\sigma^c}T^{1*}(\widetilde{K^{*}}^{-1}\widetilde{K^{*}}(f))-T^2_{\sigma^c}T^{1*}(\widetilde{K^{*}}^{-1}\widetilde{K^{*}}(f))\| \notag\\
& \quad  -\|T^2_{\sigma^c}T^{1*}(\widetilde{K^{*}}^{-1}\widetilde{K^{*}}(f))-T^2_{\sigma^c}T^{2*}(\widetilde{K^{*}}^{-1}\widetilde{K^{*}}(f))\| \notag\\
&\geq\sqrt{\alpha A_1}\|K^{*}(f)\|-\|T^{1*}\|\|T^1_{\sigma^c}-T^2_{\sigma^c}\|\|\widetilde{K^{*}}^{-1}\|\|\widetilde{K^{*}}(f)\| \notag\\
& \quad -\|T^2_{\sigma^c}\| \|T^{1*}-T^{2*}\|\|\widetilde{K^{*}}^{-1}\|\|\widetilde{K^{*}}(f)\| \notag\\
&\geq\sqrt{\alpha A_1}\|K^{*}(f)\|-\sqrt{B_1}\left(\lambda+\mu\sqrt{B_1}+\nu\sqrt{B_2}\right)\|\widetilde{K^{*}}^{-1}\|\|K^{*}(f)\| \notag\\
& \quad  -\sqrt{B_2}\left(\lambda+\mu\sqrt{B_1}+\nu\sqrt{B_2}\right)\|\widetilde{K^{*}}^{-1}\|\|K^{*}(f)\| \notag\\
&= \Big [\sqrt{\alpha A_1} - \Big(\sqrt{B_1}+\sqrt{B_2}\Big)\left(\lambda+\mu\sqrt{B_1}+\nu\sqrt{B_2}\right)\|\widetilde{K^{*}}^{-1}\|\Big] \|K^{*}(f)\|.
\end{align}
Now,  from  Proposition \ref{pro1}, the family   $\{\phi_{1j}\}_{j\in\sigma}\bigcup\{\phi_{2j}\}_{j\in\sigma^c}$ is a Bessel sequence in $\mathcal{H}$ with universal Bessel bound $B_1+B_2$. Let $T$ be its associated pre-frame operator. Then, by using \eqref{eq2.11}, we have
\begin{align*}
&\Big [\sqrt{\alpha A_1} - \Big(\sqrt{B_1}+\sqrt{B_2}\Big)\left(\lambda+\mu\sqrt{B_1}+\nu\sqrt{B_2}\right)\|\widetilde{K^{*}}^{-1}\|\Big]^2\|K^{*}(f)\|^2\\
&\leq\Big\|\sum_{j\in \sigma}\langle f, \phi_{1j}\rangle\phi_{1j}+\sum_{j\in \sigma^c}\langle f, \phi_{2j}\rangle\phi_{2j}\Big\|^2\\
&=\|TT^*f\|^2\\
&\leq\|T\|^2\|T^*f\|^2\\
&\leq(B_1+B_2)\Big(\sum_{j\in \sigma}|\langle f, \phi_{1j}\rangle|^2+\sum_{j\in \sigma^c}|\langle f, \phi_{2j}\rangle|^2\Big) \ \text{for all} \ f \in \mathcal{H}.
\end{align*}
Hence  $\Phi_1$  and $\Phi_2$ are $K$-woven with the required universal  bounds.
\endproof

\bibliographystyle{amsplain}

\end{document}